\documentclass{mpi2015-cscpreprint}

\usepackage[american]{babel}
\usepackage{amssymb}
\usepackage{amsthm,amsmath}
\usepackage{graphicx}

\newtheorem{proposition}{Proposition}

\numberwithin{equation}{section}

\begin{document}

\title{Dynamic Modelling of a Controlled Orthotropic Plate: Analytic and Data-Driven Approaches in the Frequency Domain}

\author[2]{Alexander Zuyev$^{1,3}$, Francesco Pellicano$^2$, Antonio Zippo$^2$, Giovanni~Iarriccio}
\affil[1]{Max Planck Institute for Dynamics of Complex Technical Systems, Magdeburg, Germany}
\affil[2]{Department of Engineering ``Enzo Ferrari'', University of Modena and Reggio Emilia, Modena, Italy}
\affil[3]{Institute of Applied Mathematics and Mechanics, National Academy of Sciences of Ukraine}



\abstract{
This paper is devoted to the mathematical modelling of a vibrating orthotropic plate equipped with a laminated piezosensor, under the influence of a lumped force actuation. We employ the Kirchhoff plate theory to derive the corresponding partial differential equation, assuming free boundary conditions. Analytical solutions for this boundary value problem are explored in the form of series expansions, using products of Krylov functions. Utilizing Galerkin’s method, this mathematical model is transformed into an infinite-dimensional control system characterized by modal coordinates. The transfer function of such a system is explicitly evaluated in the single-input single-output case. The computation of coefficients for finite-dimensional approximate systems is formalized in an algorithm with an arbitrary number of degrees of freedom.
Our numerical study confirms that the modeled input-output behavior shows acceptable agreement over the given frequency range.
}

 \novelty{The key contributions of this paper are as follows:
\begin{itemize}
\item mathematical model of a vibrating orthotropic plate with non-collocated actuator and sensor,
described by the Kirchhoff partial differential equation with free or simply supported boundary conditions;
\item vectorized form of the system of linear algebraic equations for approximating
solutions to the spectral problem based on the Krylov functions expansion;
\item analytical derivation of the transfer function for Galerkin's approximation of the single-input single-output orthotropic plate model;
\item evaluation of the transfer function and identification of modal frequencies using experimental data.
\end{itemize}}

\maketitle
\section{Introduction}\label{sec_intro}

Flexible structures containing composite plates and shells are widely used in modern mechanical engineering.
Analyzing the vibrations in such structures requires accurate modeling that accounts for anisotropic material properties.
While finite element methods (FEM) offer good accuracy for the modal analysis of anisotropic plates, the numerical results typically do not provide sufficient insight into the influence of mechanical parameters and the effects of time-varying input actions.
These effects are essential for addressing control-theoretic problems, particularly in relation to the task of active vibration suppression within dynamic environments.
This paper aims at comparing results of analytic and data-driven modeling of a controlled composite plate, which is controlled by a point force actuation.

The rest of this paper is organized as follows. A mathematical model of a controlled orthotropic Kirchhoff plate with free and simply supported edges is proposed in Section~\ref{sec_Kirchhoff}.
The spectral problem of the associated fourth-order differential operator is studied in Section~\ref{sec_spectral},
and the transfer function is derived in Section~\ref{sec_tf} with the use of Galerkin's method.
We justify the obtained analytical results with experimental data in Section~\ref{sec_simulations}.
The input-output behavior of a finite-dimensional models is compared with the results from the data-driven identification of the transfer function, conducted for the experimental setup at the Department of Engineering ``Enzo Ferrari'' of the University of Modena and Reggio Emilia.

\section{Controlled Kirchhoff plate model}\label{sec_Kirchhoff}
The vibration of a rectangular orthotropic plate of uniform thickness $h$ with a shaker is modeled by
the partial differential equation
\begin{equation}\label{Kirchhoff}
\begin{aligned}
\frac{\partial^2 w}{\partial t^2} &+ \alpha \frac{\partial w}{\partial t} + d_{11} \frac{\partial^4 w}{\partial x_1^4}    +d_{22} \frac{\partial^4 w}{\partial x_2^4} \\
&+ 2(d_{12}+2 d_{66})\frac{\partial^4 w}{\partial x_1^2 \partial x_2^2} =  \frac{F}{\rho h} \delta_{S_0}(x),
\end{aligned}
\end{equation}
where the function $w=w(t,x)$, $t\ge 0$, $x=(x_1,x_2)\in D = [0,\ell_1]\times [0,\ell_2]$ describes the plate middle surface (so that $w=0$ is the reference undeformed configuration).
For the sake of brevity in notation,  we will also use the subscripts in $w_t$ and $w_{x_i}$ for derivatives with respect to $t$ and $x_i$, respectively.
Equation~\eqref{Kirchhoff} is derived within the framework of Kirchhoff plate theory~\cite{amb,lagnese,TY15,MKA20,zuyev2015partial}, taking into account the force $F$ generated by the shaker.

The mechanical parameters in~\eqref{Kirchhoff} are defined as follows:
\begin{equation}\label{mechparams}
\begin{aligned}
& d_{11} = \frac{E_1 h^2}{12\rho(1-\nu_1 \nu_2)}, \; d_{22}= \frac{E_2 h^2}{12\rho(1-\nu_1 \nu_2)}, \\\
& d_{12} = \nu_2 d_{11}=\nu_1 d_{22},\;  d_{66} = \frac{G h^2}{12 \rho},
\end{aligned}
\end{equation}
where $\rho$ is the density (mass per unit volume),
$E_1$ and $E_2$ are Young's moduli, $G$ is the shear modulus,
$\nu_1$ and $\nu_2$ are Poisson's ratios such that $\nu_1 E_2 = \nu_2 E_1$,
and $\alpha$ is the viscous damping coefficient~\cite{S55}.
The considered model does not take into account the membrane forces~\cite{amb} and structural damping effects~\cite{LLD94}.

The force~$F$ acts on the shaker's contact zone $S_0\subset D$, so that the function $\delta_{S_0}(x_1,x_2)\ge 0$ vanishes outside $S_0$ and $\int_{S_0}\delta_{S_0}(x_1,x_2)  dx_1 dx_2 = 1$. If $S_0=\{(x_1^0,x_2^0)\}$ is a singleton, then $\delta_{S_0}(x_1,x_2)$ is the Dirac delta concentrated in $(x_1^0,x_2^0)$, and solutions of the Kirchhoff plate equation are treated in the generalized sense.

If the plate is simply supported at $x_1=0$ and $x_1=\ell_1$, then the corresponding boundary conditions are:
\begin{equation}\label{BC_ss}
w=0,\; w_{x_1 x_1} + \nu_2 w_{x_2 x_2} = 0 \quad\text{at}\; x_1\in\{0,\ell_1\}.
\end{equation}
In the case of  free edges $x_1=0$ and $x_1=\ell_1$, we have~\cite{BL01,lagnese}:
\begin{equation}\label{BC_F}
\begin{aligned}
& w_{x_1 x_1} + \nu_2 w_{x_2 x_2} = 0, \\
& w_{x_1 x_1 x_1}+\left( \nu_2 + \frac{4(1-\nu_1\nu_2)G}{E_1} \right) w_{x_1 x_2 x_2} =0\;
\text{at}\;\; x_1\in\{0,\ell_1\}.
\end{aligned}
\end{equation}
Similarly, the boundary conditions at free edges $x_2=0$ and $x_2=\ell_2$ are:
\begin{equation}\label{BC_x2}
\begin{aligned}
&\nu_1 w_{x_1 x_1} + w_{x_2 x_2} = 0, \\
& \left( \nu_1 + \frac{4(1-\nu_1\nu_2)G}{E_2} \right) w_{x_1 x_1 x_2} + w_{x_2 x_2 x_2} =0 \;\;\text{at}\; x_2\in\{0, \ell_2\}.
\end{aligned}
\end{equation}
The Kirchhoff plate equation~\eqref{Kirchhoff} together with the above boundary conditions can be derived from Hamilton's principle.
Indeed, let the functions $w\in C^4([t_1,t_2]\times D)$ and $F\in C([t_1,t_2])$ satisfy the following variational formulation:
 \begin{equation}\label{varform}
 \begin{aligned}
&\int_{t_1}^{t_2}\int_D \Bigl( (w_{tt} + \alpha w_t)\delta w + d_{11} w_{x_1x_1} \delta w_{x_1x_1} + d_{22} w_{x_2x_2} \delta w_{x_2x_2}\\
&+ d_{12}(w_{x_1x_1} \delta w_{x_2x_2} + w_{x_2x_2} \delta w_{x_1x_1}) + 4d_{66} w_{x_1x_2} \delta w_{x_1x_2}\Bigr) dx dt \\
&= \frac{1}{\rho h}\int_{t_1}^{t_2} \int_D F(t) \delta_{S_0}(x) \delta w(t,x)dx dt
\end{aligned}
\end{equation}
for all admissible variations $\delta w\in C^2([t_1,t_2]\times D)$ such that $w(t_1,x)=w(t_2,x)=0$ for $x\in D$.
Then, by applying integration by parts to equation~\eqref{varform} and referring to the fundamental lemma of the calculus of variations,
we observe that $w(t,x)$ satisfies differential equation~\eqref{Kirchhoff}, along with the boundary conditions~\eqref{BC_F} and~\eqref{BC_x2}.
The variational approach also allows to study general anisotropic plate models~\cite{berd} (cf.~\cite{ZS15}).
The simply supported plate at $x_1=0$ and $x_1=\ell_1$ with free edges $x_2=0$ and $x_2=\ell_2$ was considered in~\cite{ZY21}.

A strain gauge is attached to the plate surface $(x_1,x_2,w(t,x_1,x_2)+h/2)$ with $(x_1,x_2)\in C_0\subset D$ to measure deformation in the $x_1$-direction.
So, we assume that the output of the considered control system is
$$
y(t) = \int_{C_0} w_{x_1 x_1} (t,x) dx.
$$
If $C_0=\{\bar x\}$ is a singleton, then we assume that
\begin{equation}\label{output}
y(t) = w_{x_1 x_1} (t,\bar x).
\end{equation}
From the physical viewpoint, the value of $y(t)$ is proportional to the output voltage from the strain gauge.

\section{Separation of variables and spectral problem}\label{sec_spectral}
In order to study spectral properties of the considered mathematical model, we first substitute $w(t,x)=q(t) W(x)$ into~\eqref{Kirchhoff},~\eqref{BC_F},~\eqref{BC_x2} with $F=0$.
As a result, we obtain the ordinary differential equation
\begin{equation}
\ddot q(t) + \alpha \dot q(t) = -\lambda q(t)
\label{q_ODE}
\end{equation}
together with the following boundary value problem for $W(x)$:
\begin{equation}\label{BVP}
({\cal L} W)(x) = \lambda W(x),\quad W\in D({\cal L})\subset L^2(D).
\end{equation}
Here, $L^2(D)$ denotes the complex Hilbert space of square-integrable functions equipped with the standard inner product $\left<W^i,W^j\right>_{L^2(D)} = \int_D W^i(x)\overline{W^j(x)} dx$ for $W^1,W^2\in L^2(D)$, and the linear differential operator ${\cal L}:D({\cal L})\to L^2(D)$ is defined by the rule
\begin{equation}\label{L_op}
{\cal L} W =  d_{11}W_{x_1 x_1 x_1 x_1} + 2(d_{12}+2 d_{66}) W_{x_1 x_1 x_2 x_2} + d_{22} W_{x_2 x_2 x_2 x_2},
\end{equation}
where
$$
D({\cal L})= \left\{W\in H^4(D)\,:\, \begin{aligned}& W_{x_1 x_1} + \nu_2 W_{x_2 x_2} = 0\;\text{at}\;x_1\in\{0,\ell_1\},\\ &W_{x_1 x_1 x_1}+\left( \nu_2 + \frac{4d_{66}}{d_{11}} \right) W_{x_1 x_2 x_2} =0\;\text{at}\;x_1\in\{0,\ell_1\},\\
&\nu_1 W_{x_1 x_1} + W_{x_2 x_2} = 0\;\text{at}\;x_2\in\{0,\ell_2\},\\  &\left( \nu_1 + \frac{4d_{66}}{d_{22}} \right) W_{x_1 x_1 x_2} + W_{x_2 x_2 x_2} =0\;\text{at}\;x_2\in\{0,\ell_2\}\end{aligned}\right\},
$$
and $H^k(D)$ is the Sobolev space of functions $f\in L^2(D)$ such that, for each $0\le \alpha_1+\alpha_2 \le k$, the weak derivative $\frac{\partial^{\alpha_1+\alpha_2}f}{\partial^{\alpha_1}x_1 \partial^{\alpha_2}x_2}$ exists and belongs to $L^2(D)$.
We introduce the symmetric bilinear form ${\cal B}: H^2(D)\times H^2(D)\to \mathbb C$ such that
\begin{equation}\label{B_form}
\begin{aligned}
{\cal B}(W^i,& W^j) = \int_D
\Bigl(d_{11} W^i_{x_1x_1} W^j_{x_1x_1} + d_{22} W^i_{x_2x_2} W^j_{x_2x_2}\\
&+ d_{12}(W^i_{x_1x_1} W^j_{x_2x_2} + W^i_{x_2x_2} W^j_{x_1x_1}) + 4d_{66} W^i_{x_1x_2} W^j_{x_1x_2}\Bigr) dx.
\end{aligned}
\end{equation}
By exploiting the integration by parts in~\eqref{B_form}, it is easy to see that
\begin{equation}\label{B_diffrelation}
{\cal B}(W^i,W^j) = \int_D ({\cal L} W^i)(x) W^j(x)dx,
\end{equation}
provided that $ W^i\in D({\cal L})$ and $W^j\in H^2(D)$.
Some important properties of eigenvalues and eigenfunctions are summarized below.
\begin{proposition}\label{lem_orthogonality}
The solutions to boundary value problem~\eqref{BVP} satisfy the following properties.
\begin{enumerate}
\item[1)] All eigenvalues $\lambda$ of boundary value problem~\eqref{BVP} are real.
\item[2)] Let $(\lambda_i,W^i)$ and $(\lambda_j,W^j)$ be solutions of problem~\eqref{BVP}.
If $\lambda_i\neq \lambda_j$, then $\left<W^i,W^j\right>$ = 0.
\end{enumerate}
\end{proposition}

{\textbf Proof.}
Let $W(x)\not\equiv 0$ be a solution of boundary value problem~\eqref{BVP} with some $\lambda\in\mathbb C$. Then, because of~\eqref{B_diffrelation}, we conclude that
$$
\begin{aligned}
\lambda \left<W,W\right>_{L^2(D)} &= \left<{\cal L} W,W\right>_{L^2(D)} = {\cal B}(W,\overline{W}) = {\cal B}(\overline{W},W) \\
&= \overline{{\cal B}(W,\overline{W})}
 = \overline{\left<{\cal L} W,W\right>}_{L^2(D)} \\
 &= \overline{\lambda  \left<W,W\right>}_{L^2(D)} = \overline{\lambda}  \left<W,W\right>_{L^2(D)}.
\end{aligned}
$$
Hence, $\lambda=\overline{\lambda}$, which proves assertion~1).

To prove assertion~2) for $W^i,W^j\in D({\cal L})$ such that ${\cal L} W^i = \lambda_i W^i$ and ${\cal L} W^j = \lambda_j W^j$, we employ identity~\eqref{B_diffrelation} and the equalities
$$
\begin{aligned}
\lambda_i \left<W^i,W^j\right>_{L^2(D)} &=  {\cal B} (W^i,\overline{W^j}) = \overline{{\cal B} (W^j,\overline{W^i})} \\
&= \overline{\int_D ({\cal L} W^j)(x) \overline{W^i}(x)dx}
 = \overline{\lambda_j \left<W^j,W^i\right>}_{L^2(D)} \\ & = \overline{\lambda_j} \left<W^i,W^j\right>_{L^2(D)} = {\lambda_j} \left<W^i,W^j\right>_{L^2(D)}.
\end{aligned}
$$
Therefore, $(\lambda_i - \lambda_j)\left<W^i,W^j\right>_{L^2(D)} =0$, so that the eigenfunctions $W^i$ and $W^j$ are orthogonal in $L^2(D)$ whenever $\lambda_i\neq \lambda_j$.
$\square$

The boundary value problem, as stated in~\eqref{BVP}, can be reformulated into its weak form by employing identity~\eqref{B_diffrelation}.
This leads to the following formulation: find a $\lambda\in \mathbb R$
and $W\in D({\cal L})$ such that
\begin{equation}\label{BVP_weak}
{\cal B}(W,V) = \lambda \int_D W(x) V(x) dx,
\end{equation}
for each $V\in H^2(D)$.

We will search for solutions of~\eqref{BVP_weak} in the form of series
\begin{equation}
W(x) = \sum_{n,m=0}^\infty c_{nm} \phi_n\left(\frac{x_1}{\ell_1}\right) \phi_m\left(\frac{x_2}{\ell_2}\right),\quad x=(x_1,x_2)\in D,
\label{W_series}
\end{equation}
where $\{\phi_n(x)\}_{n=0}^\infty$ are normalized eigenfunctions of the Euler--Bernoulli beam model with free ends.
Namely, the functions $\phi_n(x)$ are solutions of the boundary value problem
\begin{equation}
\begin{aligned}
& \phi_n''''(x) = \beta_n^4 \phi_n(x),\quad x\in (0,1),\\
& \phi_n''(0)=\phi_n''(1)=0,\\
& \phi_n'''(0)=\phi_n'''(1)=0.
\end{aligned}
\label{BVP_EB}
\end{equation}
The eigenvalues $\beta_k^4$ are such that $\beta_0=\beta_1=0$, and $\beta_n>0$ satisfy
$$
\cos \beta_n \cdot \cosh \beta_n = 1\quad \text{for}\; n=2,3,... \; .
$$
The above transcendent equation implies that
$$
\beta_n = \frac{\pi (2n-1)}{2} + o \left( \frac1n\right) \quad \text{as}\; n \to \infty.
$$

The normalized eigenfunctions of~\eqref{BVP_EB} are represented in the form
\begin{equation}
\begin{aligned}
\phi_0(x) &= 1,\; \phi_1(x) = \sqrt{3}(2x-1),\\
\phi_n(x) &= \frac{\psi_n(x)}{\|\psi_n(x)\|_{L^2(0,1)}},\\
\psi_n(x) &= \sin \beta_n x + \sinh \beta_n x - \frac{\cosh \beta_n - \cos \beta_n}{\sinh \beta_n + \sin \beta_n}\\
&
\times\left( \cos\beta_n x + \cosh \beta_n x\right)\;\; \text{for}\; n=2,3,...,
\end{aligned}
\label{phi_n}
\end{equation}
such that $\|\phi_n\|_{L^2(0,1)}=1$
and $\left<\phi_n,\phi_m\right>_{L^2(0,1)}=0$ for all $n\neq m$, where $\left<\cdot,\cdot\right>_{L^2(0,1)}$ stands for the standard inner product in the Hilbert space $L^2(0,1)$.

To determine the coefficients $c_{nm}$ in~\eqref{W_series}, we apply Galerkin's method.
This requires that equation~\eqref{BVP_weak} is satisfied for each function $V$ of the form
$$
\phi_{kj}(x) = \phi_k\left(\frac{x_1}{\ell_1}\right) \, \phi_j\left(\frac{x_2}{\ell_2}\right).
$$
In other words, we assume that
\begin{equation}\label{BVP_Galerkin}
{\cal B}(W,\phi_{kj}) = \lambda \int_D W(x) \phi_{kj}(x) dx\quad\text{for}\;\;k,j=0,1,2,\, ...\, .
\end{equation}
After substituting the expansion given in~\eqref{W_series} into equation~\eqref{BVP_Galerkin} and utilizing the orthogonality of the functions $\phi_j$, we formally derive
\begin{equation}
\begin{aligned}
&\left( \frac{d_{11} \beta_k^4}{\ell_1^4} + \frac{d_{22} \beta_j^4}{\ell_2^4} - \lambda \right)c_{kj} + \frac{d_{12}}{\ell_1^2 \ell_2^2} \sum_{n,m=0}^\infty (\varkappa_{nk} \varkappa_{jm} + \varkappa_{kn} \varkappa_{mj}) c_{nm}\\
& + \frac{4 d_{66}}{\ell_1^2 \ell_2^2} \sum_{n,m=0}^\infty \theta_{nk} \theta_{mj}c_{nm}= 0,
\end{aligned}
\label{algsys}
\end{equation}
where
\begin{equation}
\varkappa_{nk}= \left<\phi_n'',\phi_k\right>_{L^2(0,1)},\; \theta_{nk} = \left<\phi_n',\phi_k'\right>_{L^2(0,1)}.
\label{kappa_def}
\end{equation}
Note that $\varkappa_{nj}=0$ for $n<2$, and $\theta_{0k}=\theta_{n0}=0$ due to formulas~\eqref{phi_n}.

The infinite system of linear algebraic equations~\eqref{algsys} with respect to $c_{kj}$ is written in the matrix form as
\begin{equation}
\frac{d_{11}}{\ell_1^4} B^4 C + \frac{d_{22}}{\ell_2^4} CB^4 +\frac{d_{12}}{\ell_1^2 \ell_2^2}\left(\varkappa C \varkappa + \varkappa^T C^T \varkappa^T\right)+  \frac{4 d_{66}}{\ell_1^2 \ell_2^2} \theta C \theta = \lambda C
\label{algsys_abstract}
\end{equation}
with the infinite matrices $B=\mathrm{diag}(\beta_0,\beta_1,\beta_2, ...)$, $C=(c_{kj})_{k,j=0}^\infty$,
$\varkappa= (\varkappa_{kj})_{k,j=0}^\infty$, $\theta= (\theta_{kj})_{k,j=0}^\infty$.
A finite-dimensional approximation of this system can be obtained by fixing a natural number $N$ and truncating
the entries and summands with indices $k,j>M$ in~\eqref{algsys_abstract}. As a result, we obtain the following truncated version of system~\eqref{algsys_abstract}:
\begin{equation}
\begin{aligned}
\frac{d_{11}}{\ell_1^4} B_N^4 \tilde C &+ \frac{d_{22}}{\ell_2^4} \tilde C B_N^4 +\frac{d_{12}}{\ell_1^2 \ell_2^2}\left(\varkappa_N \tilde C \varkappa_N + \varkappa_N^T {\tilde C}^T \varkappa_N^T\right)\\
&+  \frac{4 d_{66}}{\ell_1^2 \ell_2^2} \theta_N {\tilde C} \theta_N = \lambda \tilde C,
\end{aligned}
\label{algsys_truncated}
\end{equation}
where $B_N=\mathrm{diag}(\beta_0,\beta_1, ...,\beta_N)$, $\tilde C=(\tilde c_{k,j})_{k,j=0}^N$, $\varkappa_N=(\varkappa_{k,j})_{k,j=0}^N$, $\theta_N= (\theta_{kj})_{k,j=0}^N$  are square matrices of size $(N+1)\times (N+1)$.
The vectorized form of matrix equation~\eqref{algsys_truncated} is
\begin{equation}
M \hat C = \lambda \hat C,
\label{sys_vectorized}
\end{equation}
where
\begin{equation}
\begin{aligned}
M = &\frac{d_{11}}{\ell_1^4} (I \otimes B_N^4) + \frac{d_{22}}{\ell_2^4} (B_N^4 \otimes I) + \frac{4d_{66}}{\ell_1^2 \ell_2^2} ( \theta_N^T \otimes \theta_N)\\
&+ \frac{d_{12}}{\ell_1^2 \ell_2^2}(\varkappa_N^T \otimes \varkappa_N+ (\varkappa_N \otimes \varkappa_N^T)K)
\end{aligned}
\label{M_matrix}
\end{equation}
is the square matrix of size $(N+1)^2\times (N+1)^2$, the column vector $\hat C = \mathrm{vec}\, \tilde C$ is the vectorization of $\tilde C$,
$K$ is the commutation matrix such that $K \mathrm{vec}\, \tilde C = \mathrm{vec}\, (\tilde C^T)$,
and ``$\otimes$'' stands for the Kronecker product.
Thus, we will treat the eigenvalues $\lambda$ of the matrix $M$ given by~\eqref{M_matrix} as approximations to the eigenvalues of problem~\eqref{BVP},
and the corresponding truncated expansions
\begin{equation}\label{W_truncated}
\tilde W(x) =
\sum_{k,j=0}^N \tilde c_{kj} \phi_k\left(\frac{x_1}{\ell_1}\right) \phi_j\left(\frac{x_2}{\ell_2}\right)
\end{equation}
will be considered as approximate solutions of~\eqref{BVP}.

\section{Transfer function}\label{sec_tf}
In order to study the input-output behavior of the control system~\eqref{Kirchhoff},~\eqref{BC_F},~\eqref{BC_x2},~\eqref{output}, we evaluate the transfer function for its Galerkin approximation.
 Specifically, we assume that a finite set of nonnegative integers $\cal N$ describes the indices of relevant modes of vibrations and denote by
$(\lambda_n,W_n)$ the nontrivial solutions of the corresponding spectral problem~\eqref{BVP} for $n\in\cal N$.
Then, assuming that all $\lambda_n$ are distinct and using Proposition~1, we deduce that the transfer function for Galerkin's approximation of~\eqref{Kirchhoff},~\eqref{BC_F},~\eqref{BC_x2},~\eqref{output} takes the following form:
\begin{equation}
H_{\cal N}(s) = \frac{1}{\rho h}\sum_{n\in\cal{N}} \frac{W_{n_{x_1 x_1}}(C_0) W_n(S_0)}{(s^2+\alpha s + \lambda_n)\|W_n\|^2_{L^2(D)}}.
\label{TF_analytic}
\end{equation}
A numerical approximation of $H_{\cal N}(s)$ can be derived by substituting the eigenvalues of $M$ and the corresponding truncated expansions $\tilde W(x)$ from~\eqref{W_truncated} into~\eqref{TF_analytic}.

\section{Experimental results and numerical simulations}\label{sec_simulations}
We will present experimental results for a three-layer composite plate  consisting of two graphite epoxy sheets bonded through a NOMEX honeycomb layer. The carbon fiber stacking sequence is 0°/90°.
We represent this plate by an orthotropic Kirchhoff model with free boundary conditions and the following mechanical parameters:
$$
\begin{aligned}
& \ell_1 = 1\, m,\; \ell_2= 0.5\, m,\, h=3.6\,mm,\, \rho = 505.6\,\frac{kg}{m^3},\\
& E_1 = 23\,GPa,\;E_2 = 14\,GPa,\, G=2.2\,GPa,\\
& \nu_1 = 0.25,\,\nu_2 = \frac{\nu_1 E_2}{E_1}.
\end{aligned}
$$
The shaker is attached at point $S_0$ with coordinates $(0.17 \,m,0.25 \, m)$, and we assume that the piezosensor at point $C_0$ with coordinates $(0.5 \,m,0.21\,m)$ provides measurements of $w_{x_1x_1}$.

The force actuation at $S_0$ is measured by a load cell, and the signals are processed by a dSpace control/acquisition system.
In this experimental setup, we use pulse actuation of the shaker to record the input signal of the plate (force applied by the shaker at $S_0$, measured by the load cell) and the output (voltage from the piezosensor at $C_0$).
The resulting data are recorded into Matlab files with a sampling rate of $30\, kHz$.

To compare theoretical modal frequencies of the orthotropic Kirchhoff plate model with experimental data,
we apply Galerkin's method described in Section~\ref{sec_spectral} to derive truncated expansions of the modal shapes given by~\eqref{W_truncated} for $N=6$.
The matrices that define the matrix $M$ in~\eqref{M_matrix} are as follows:
{
$$
B_N = \mathrm{diag} (0,0,4.730 7.853, 10.996, 14.137, 17.279),
$$
$$
\varkappa_N =\begin{pmatrix}
0 & 0 & 0 & 0 & 0 & 0 & 0 \\
0 & 0 & 0 & 0 & 0 & 0 & 0 \\
-18.59 & 0 & -12.30 & 0 & 1.80 & 0 & 0.57 \\
0 & 40.59 & 0 & -46.05 & 0 & 5.29 & 0 \\
-43.98 & 0 & 52.58 & 0 & -98.91 & 0 & 9.86 \\
0 & 84.09 & 0 & 55.51 & 0 & -171.59 & 0 \\
-69.12 & 0 & 101.62 & 0 & 60.13 & 0 & -263.99
\end{pmatrix},
$$
$$
\theta_N = \begin{pmatrix}
0 & 0 & 0 & 0 & 0 & 0 & 0 \\
0 & 12 & 0 & 13.86 & 0 & 13.864 & 0 \\
0 & 0 & 49.48 & 0 & 0 & 35.38 & 0 & 36.61 \\
0 & 13.86 & 0 & 108.93 & 0 & 57.59 & 0 \\
0 & 0 & 35.38 & 0 & 186.87 & 0 & 78.10 \\
0 & 13.86 & 0 & 57.59 & 0 & 284.68 & 0 \\
0 & 0 & 36.61 & 0 & 78.10 & 0 & 402.23
\end{pmatrix}.
$$
}

If $\lambda_n$ is an eigenvalue of the matrix $M$ in~\eqref{M_matrix},
then the corresponding modal frequency in Hertz is $\nu_n=\frac{\sqrt{\lambda_n}}{2\pi}$.
The lowest five modal frequencies, corresponding to nonzero eigenvalues of $M$, are presented in Table~1 alongside experimentally identified frequencies of the plate.

\begin{table}[h!]
\begin{center}
\begin{tabular}{|c|c|c|}
\hline
$n$ & Experimental frequency, Hz &  Theoretical frequency, Hz \\ \hline
1 & 16.2 & 16.0 \\  \hline
2 & 25.1 & 25.4 \\ \hline
3 & 40.8 & 41.4 \\ \hline
4 & 68.0 & 70.1 \\ \hline
5 & 83.0 & 81.0 \\ \hline
\end{tabular}
\caption{Experimental and theoretical frequencies.}
\end{center}
\label{tab1}
\end{table}

We observe acceptable agreement between the theoretically and experimentally evaluated frequencies of the considered plate setup. It should be noted that the theoretical values in Table~1 are computed without taking into account the damping parameter $\alpha$, which appears in the equation for modal amplitudes~\eqref{q_ODE} together with $\lambda$.
Thus, further refinement of these estimates is possible upon identifying the damping parameter.

\begin{figure}[h!]
 \includegraphics[width=1\linewidth] {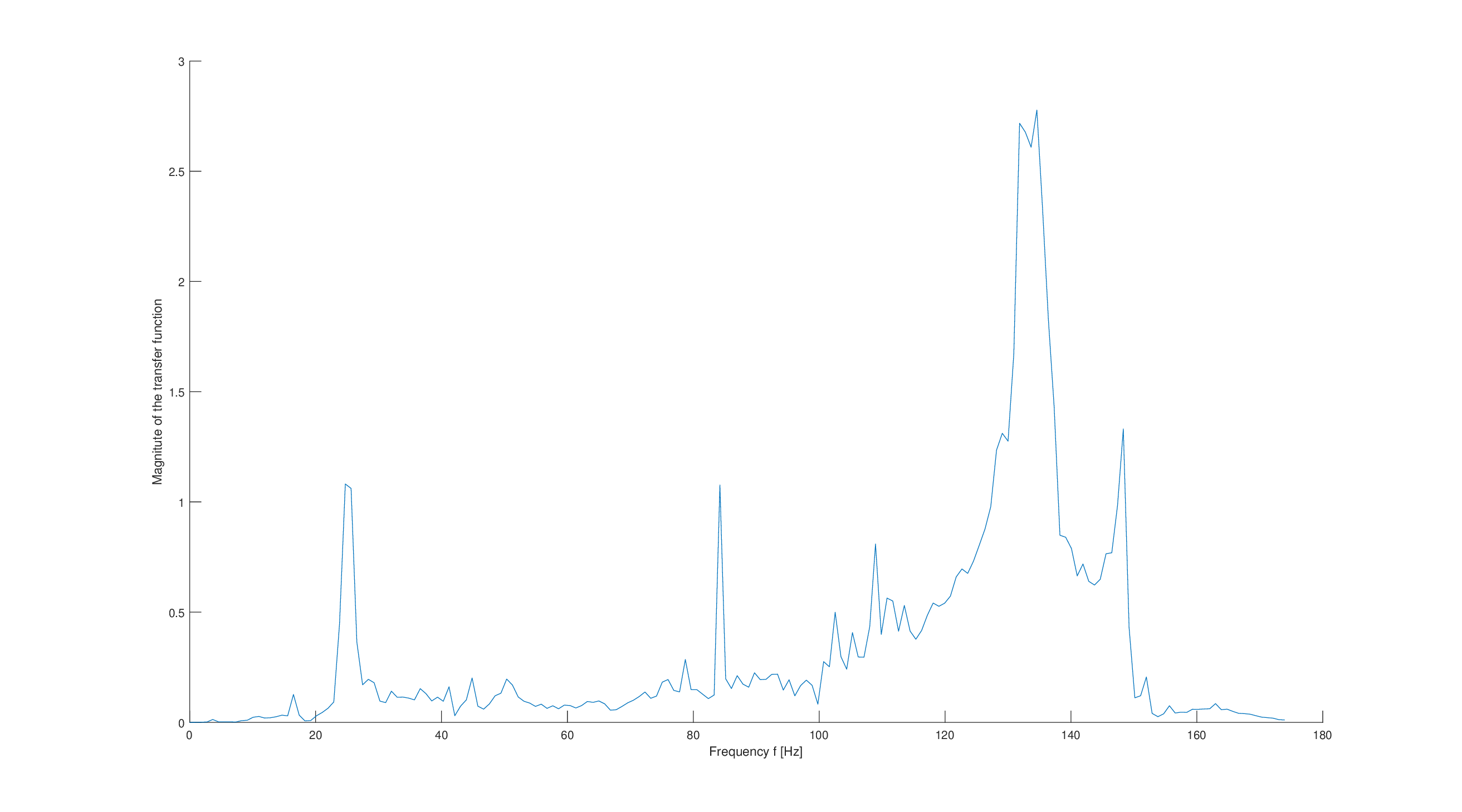}
 \caption{Experimental transfer function.}
 \label{fig1}
\end{figure}

\begin{figure}[h!]
 \includegraphics[width=1\linewidth] {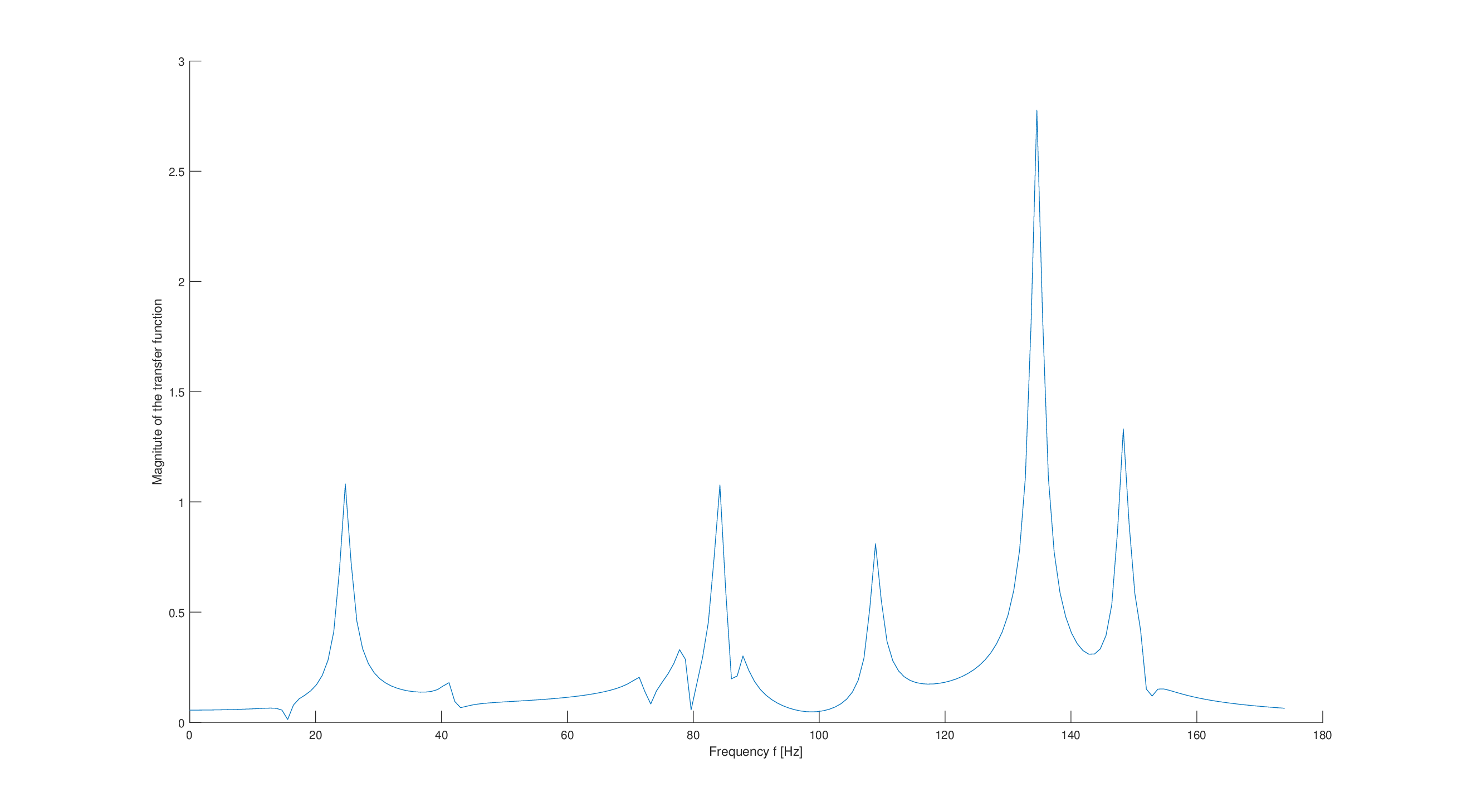}
 \caption{Theoretical transfer function for Galerkin's approximation.}
 \label{fig2}
\end{figure}

We apply the fast Fourier transform to the measured input and output signals, using with $2^{15}$ sample points, to evaluate the transfer function from the experimental data. The resulting transfer function is plotted in Fig.~\ref{fig1}.
We also find eigenvalues of the matrix $M$ in~\eqref{sys_vectorized} to evaluate the model-based transfer function.
The matrix $M$ has three zero eigenvalues $\lambda_0=\lambda_1=\lambda_2=0$ due to the free edged of the considered plate model.
By considering the set of nonzero eigenvalues $\lambda_3,\lambda_4,..., \lambda_{12}$ of $M$ for ${\cal N}=\{3,4,...,12\}$ and computing the corresponding approximate modal shapes using~\eqref{W_truncated},
we evaluate the model-based transfer function $H_{\cal N}$ of the considered plate model by formula~\eqref{TF_analytic}.
This transfer function is depicted in Fig.~\ref{fig2}, showing the frequency range from $0$ to $175$ Hz.

\newpage
\section{Conclusion}\label{sec_conclusion}
We have derived the mathematical model for an orthotropic vibrating plate, actuated by a point force and equipped with a piezoelectric sensor.
The boundary value problem, described in Section~\ref{sec_Kirchhoff}, accommodates plates with both simply supported and free edges.
For a plate with free edges, the displacement representation~\eqref{W_series}, using Krylov functions, provides a constructive approach
to solving the associated  spectral problem through the infinite matrix equation~\eqref{algsys_abstract}.
A finite-dimensional approximation of this system is effectively implemented in a vectorized form (equation~\eqref{sys_vectorized}) and has been successfully tested using Matlab.
The validity of our analytical model and computational approach is further supported by the experimental results presented in Section~\ref{sec_simulations}.
Future work will focus on developing data-driven methods for inferring reduced-order models of plates. This includes exploring the Loewner framework~\cite{antoulas2017tutorial,zuyev2023approximating}, both in theoretical and experimental contexts with our setup.

\end{document}